\documentclass[bjps, preprint]{imsart}

\usepackage[utf8]{inputenc}
\usepackage[english]{babel}
 
\usepackage{tikz-cd}
\usetikzlibrary{arrows.meta,arrows}

\usepackage{enumerate}
\usepackage{amsmath}
\usepackage{amsthm}
\usepackage{amsfonts}
\usepackage{enumitem}
\usepackage{textcomp}
\usepackage{centernot}

\theoremstyle{plain}
\newtheorem*{theorem*}{Theorem}
\newtheorem{theorem}{Theorem}

\newtheorem{proposition}[theorem]{Proposition}

\newtheorem{remark}[theorem]{Remark}

\newcommand{\pr}{\mathbb{P}}

\newcommand{\Rp}{\mathbb{R}_+}

\newcommand{\X}{\mathbf{X}}

\begin{document}
\begin{frontmatter}

\title{A note on monotonicity of spatial epidemic models}

\runtitle{Monotonicity of spatial epidemics}

\begin{aug}
  \author{\fnms{Achillefs}  \snm{Tzioufas}\thanksref{a}\ead[label=e1]{tzioufas@ime.usp.br}}

  \runauthor{A. Tzioufas}

  \affiliation[a]{Instituto de Matem\' atica e Estat\' istica, Universidade de S\~ ao Paulo}

  \address{Rua do Mat\~ ao, 1010, CEP 05508-900- S\~ ao Paulo, Brasil.}

\end{aug}

\begin{abstract} 
The epidemic process on a graph is considered for which infectious contacts occur at rate which depends on whether a susceptible is infected for the first time or not. We show that the Vasershtein coupling extends if and only if secondary infections occur at rate which is greater than that of initial ones. Nonetheless we show that, with respect to the probability of occurrence of an infinite epidemic, the said proviso may be dropped regarding the totally asymmetric process in one dimension, thus settling in the affirmative this special case of the conjecture for arbitrary graphs due to [Stacey (2003), {\em Ann. Appl. Probab.} {\bf 13}, 669-690]. 
\end{abstract}

\begin{keyword}[class=MSC]
\kwd[Primary ]{60K35}
\kwd[; secondary ]{82C22}
\end{keyword}

\begin{keyword}
\kwd{ three state contact processes; stochastic domination; attractiveness; contact process; standard spatial epidemic}
\end{keyword}

\end{frontmatter}

\section{Introduction and results}
The \textit{three state contact process} infection rates $(\lambda,\mu)$ on a digraph $G \equiv G(V,E)$ is a continuous-time Markov process $\zeta_{t}$ on the configuration space $\X = \{-1,0,1\}^{V}$, i.e.\ the set of all functions from $V$ to $\{-1,0,1\}$. Transition rates for $\zeta_{t}$ are specified via the flip rates for $\zeta_{t}(u)$, which are given as follows:
\begin{equation*}\label{rates}
\begin{array}{cl}
-1 \rightarrow 1 & \mbox{ at rate } \lambda |\{\overrightarrow{vu} \in E :\zeta_{t}(v)= 1\}| \\
\mbox{ } 0 \rightarrow 1 & \mbox{ at rate } \mu |\{\overrightarrow{vu} \in E:\zeta_{t}(v)=1\} | \\ 
\mbox{ } 1 \rightarrow 0 & \mbox{ at rate } 1,
\end{array}
\end{equation*}
$t\geq0$, and where $|A|$ denotes cardinality of set $A$. Initial configuration $\eta$ and infection rates $(\lambda,\mu)$ are incorporated in the notation below in the fashion: $\zeta_{t}^{\{\eta, (\lambda,\mu)\}}$. The configuration space $\X$ is endowed with the usual component-wise partial order, given by writing $\eta_{1} \leq \eta_{2}$ whenever $\eta_{1}(x) \leq \eta_{2}(x)$, for all $x \in V$.  If $(Y_{t})_{t\geq0}$ and $(W_{t})_{t\geq0}$ are two stochastic processes on $\X$, we write $Y_{t} \geq_{st.} W_{t}$ to denote that $Y_{t}$ stochastically dominates $W_{t}$,  which is, that the two processes can be defined on a common probability space such that:  $Y_{t} \geq W_{t}$, for all $t\geq0$, almost surely.\footnote{Equivalently, if $\nu_{t}, \mu_{t}$ denote the distributions of $Y_{t}, W_{t}$ respectively, then $Y_{t} \geq_{st.} W_{t}$ if, for all  $t\geq 0$, 
\begin{equation}\label{stmon}
\int_{\X} f d\nu_{t} \geq \int_{\X} f d\mu_{t},  \mbox{ for all increasing } f \mbox{ on } \X,
\end{equation}
where $f \in C(\X)$, the space of continuous functions on $\X$ equipped with the uniform norm, is said to be increasing in the component-wise sense.}
The usual assumption that $G$ is of bounded degree, needed to assure uniqueness of the process, is adopted here. Note also that although the process is an interacting particle system, it is not a spin system as the contact process itself due to the inclusion of a third state. 
For background and general information about interacting particle systems, we refer to Liggett \cite{L85}, \cite{L99} and Durrett \cite{D88}, \cite{D95}.

The three-state contact process admits the following epidemiological interpretation, which will be in effect in the sequel. If $\zeta_{t}(x) = 1$, site $x$ is regarded as infected; if $\zeta_{t}(x) = -1$, it is regarded as susceptible and never infected, and if  $\zeta_{t}(x) = 0$, as susceptible and previously infected. Thus, transitions $-1 \rightarrow 1$,  $ 0 \rightarrow 1$, and $1 \rightarrow 0$ may be thought of as \textit{initial infections, secondary infections and recoveries} respectively. Observe that, when reinfections occur at the same rate as initial infections, i.e.\ $\lambda = \mu$, and  when reinfections are disallowed, i.e.\ $\mu=0$, the process reduces to the extensively studied \textit{contact process} and the \textit{standard spatial epidemic} process respectively. The three state contact process was introduced in the mathematics literature by Durrett and Schinazi \cite{DS00}, independently of Grassberger, Chate, and Rousseau \cite{GCR97} that introduced it in the physics one first. The process has been considered afterward under a different name in Stacey \cite{S03}, and in $\mathsection\mathsection$  2, 3 and 5 in \cite{T11}. Because interest in the process stems principally from understanding the variation induced in properties of the contact process when allowing for a different initial infection rate, focus is placed upon analysis of the process from initial configurations $\eta_{A}$ such that $\eta_{A}(x)=1$, for $x\in A$, $\eta_{A}(x)=-1$, otherwise, $A$ finite. In what follows, we write simply $\zeta_{t}^{\{A, (\lambda,\mu)\}}$ for $ \zeta_{t}^{\{\eta_{A}, (\lambda,\mu)\}}$, and also, $\zeta_{t}^{\{u, (\lambda,\mu)\}}$ instead of $\zeta_{t}^{\{\{u\}, (\lambda,\mu)\}}$.  Furthermore, we use the shorthand $\{\zeta_{t}^{\{\eta, (\lambda,\mu)\}} \mbox{ survives}\}$   to denote  $\{\forall \mbox{ }t, \zeta_{t}^{\{\eta, (\lambda,\mu)\}}(x) =1 \mbox{ for some } x \}$, which is, the event an infinite epidemic occurs. 



%

Our first result regards extending a stochastic monotonicity property of spin systems, the so-called \textit{attractiveness}, to the three-state contact process on $G$. This result improves Proposition 2.1 in Stacey \cite{S03}, where the additional assumption $\lambda' \geq \mu$ in our notation is required.

\begin{theorem}\label{PIPmon}
Let $\zeta_{t}' \equiv \zeta_{t}^{\{\eta', (\lambda', \mu')\}}$ and $\zeta_{t} \equiv \zeta_{t}^{\{\eta, (\lambda, \mu)\}}$ be the three state contact processes with initial configurations $\eta'$ and $\eta$ and infection rates $(\lambda', \mu')$ and $(\lambda, \mu)$, respectively.  We have that
\begin{equation}\label{monimp}
\lambda' \geq \lambda, \mu' \geq \mu  \mbox{ and } \mu'\geq \lambda \hspace{2mm} \Longrightarrow   \mbox{ for all } \eta'  \geq \eta: \zeta_{t}' \geq_{st.} \zeta_{t}. 
\end{equation}
\end{theorem}


We give in the next remark a compound form of Theorem \ref{PIPmon} which facilitates applications, and then explicate on the necessity of condition $\mu'\geq \lambda$ for (\ref{monimp}) to hold. 

\begin{remark}\textup{
\textit{Monotonicity in the initial configuration}:
\begin{equation}\label{rem11}
\forall (\lambda, \mu) \mbox{ s.t. } \mu \geq \lambda, \mbox{ }  \eta \leq \eta'  \hspace{2mm} \Longrightarrow \mbox{ }\zeta_{t}^{\eta'} \geq_{st.} \zeta_{t}^{\eta}. 
\end{equation}
\textit{Monotonicity in the infection rates}:
\begin{equation}\label{rem12}
\forall \eta, \mbox{ } \lambda' \geq \lambda, \mu' \geq \mu  \mbox{ and } \mu'\geq \lambda \hspace{2mm} \Longrightarrow  \mbox{ } \zeta_{t}^{(\lambda', \mu')} \geq_{st.} \zeta_{t}^{(\lambda, \mu)}. 
\end{equation}
Regarding $\zeta_{t}$ on arbitrary $G$, note that the condition $\mu \geq \lambda$ cannot be dropped for monotonicity in the initial configuration property to hold (cf.\  (\ref{rem1}) below); further, condition $\mu'\geq \lambda$ cannot be dropped for monotonicity in the infection rates property to hold  (cf.\  (\ref{rem2})). }
\end{remark}

The method of proof of Theorem \ref{PIPmon} is based on an extension of what is known in the context of spin systems as the Vasershtein  (or basic) coupling; it differs from the techniques of Stacey \cite{S03} which rely on a Harris' graphical representation, introduced in Harris \cite{HA78}, for this process. For background on stochastic monotonicity, coupling and attractiveness in particular, we refer the reader to Lindvall \cite{Lin02} and Liggett \cite{L99}, see also Durrett \cite{D81}. This important in the study of interacting particle systems property was introduced and studied in Holley \cite{H72}. For differences and similarities among the basic coupling and that yielded by Harris' graphical representations for spin systems, see Liggett \cite{L85}, Chpt.\ III, $\mathsection $6. We also remark, for ease of reference below, on a simple consequence of monotonicity in the infection rates. 
\begin{remark}\label{remcurves}\textup{Let $\theta(\lambda, \mu) = \pr(\zeta_{t} \mbox{ survives})$. The following critical curves may be defined: $\mu_{c}(\lambda) := \inf\{\mu \geq \lambda: \theta(\lambda, \mu) >0\}$ and  $\lambda_{c}(\mu) := \inf\{\lambda \leq \mu: \theta(\lambda, \mu) >0\}$.}
\end{remark}










For the next statement only we consider a specific graph; let $G_{1,+} \equiv G_{1,+}(V,E)$ be such that $V= \{0,1,\dots\}$ and $E= \{\overrightarrow{vu}: u = v+1 \}$. The contact process on $G_{1,+}$  has been extensively studied in the literature and is referred to as the one-sided basic contact process, see, for instance, Griffeath \cite{G79}, Chpt.\ II, $\mathsection$ 4, and Schonmann \cite{S86}.
Regarding the one-sided basic three state contact process, we note that (\ref{eq:5}) 
below addresses Question 5.1 in Stacey \cite{S03} and (\ref{extra}) settles in the affirmative Conjecture 5.2 there.

\begin{theorem}\label{onesided}
Let $\zeta_{t}' \equiv \zeta_{t}^{\{\eta_{o}, (\lambda', \mu')\}}$ and $\zeta_{t} \equiv \zeta_{t}^{\{\eta_{o}, (\lambda, \mu)\}}$ be the three state contact process on $G_{1, +}$ with infection rates $(\lambda', \mu')$ and $(\lambda, \mu)$ respectively, and the same initial configuration $\eta_{o}$ such that $\eta_{o}(0)=1$ and $\eta_{o}(n)=-1$, for all $n \geq1$. Let also $I'_{n} = \{t \geq0: \zeta'_{t}(n)=1\}$ and $I_{n} = \{t \geq0: \zeta_{t}(n)=1\}$. We have that
\begin{equation}\label{eq:5}
\lambda' \geq \lambda \mbox{ and } \mu' \geq \mu \mbox{ }  \Longrightarrow  | I_{n}' \cap [0,t] | \geq_{st.}  | I_{n} \cap [0,t] |. 
\end{equation}
Furthermore, 
\begin{equation}\label{extra}
\lambda' \geq \lambda \mbox{ and } \mu' \geq \mu \mbox{ }  \Longrightarrow \pr\left(\zeta_{t}' \textup{ survives}\right) \geq \pr\left( \zeta_{t} \textup{ survives}\right).
\end{equation}
\end{theorem} 

It is interesting to contrast our principle monotonicity property, given in (\ref{eq:5}) above, with that induced by the partial order above. Note first that (\ref{eq:5}) states that the total infected time before any fixed time is at least as great for $\zeta_{t}'(n)$ as it is for $\zeta_{t}(n)$. Note however that Remark \ref{conint} below yields that the following holds. 
\begin{equation}\label{conint2}
\lambda < \lambda' \mbox{ and }  \mu' = \mu = 0 \centernot\Longrightarrow \forall n: \zeta_{t}(n) \leq_{st.} \zeta_{t}'(n), 
\end{equation}
where $\centernot\Longrightarrow$ denotes that the implication is false.
Hence, the monotonicity property in (\ref{eq:5}) corresponds to a strictly weaker property to that in (\ref{conint2}). We note that this weak monotonicity property may also be relevant for other models which have the standard spatial epidemic and the contact process as special cases, see for instance, Van Den Berg, Grimmett and Schinazi \cite{BGS98}. Regarding the proof of Theorem \ref{onesided}, we note that it relies on a stochastic comparison (coupling) construction which exploits the memory-less property of the exponential distribution and the spatial restriction imposed in an essential way.
Note also that the following direct consequence of (\ref{extra}) gives an extension of the conclusion in Remark \ref{remcurves} above, which is valid however for general $G$. 

\begin{remark}
For the three state contact process on $G_{1, +} $ the following critical curves may be defined: $\mu_{c}(\lambda) := \inf\{\mu: \theta(\lambda, \mu) >0\}$ and  $\lambda_{c}(\mu) := \inf\{\lambda: \theta(\lambda, \mu) >0\}$. 
\end{remark}





A comparison between the probability of an infinite epidemic in the three state contact process and in the standard spatial epidemic, noted without proof for the case that $G$ is the $d$-dimensional lattice in Durrett and Schinazi \cite{DS00}, Proposition 2, is given in the following statement. Roughly speaking, the result provides that permitting for secondary infections cannot cause the probability of an infinite epidemic in the standard spatial epidemic process to decrease. 


\begin{proposition}\label{propff} 
Let $\zeta_{t}^{\{w,(\lambda,\mu)\}}$ and $\xi_{t}^{\{w, \lambda\}} \equiv \zeta_{t}^{\{w,(\lambda,0)\}}$ be the three state contact process parameters $(\lambda,\mu)$ and the standard spatial epidemic parameter $\lambda$ respectively, and the same initial configuration $\eta_{w}$, $w \in V$. We have that
\begin{equation}\label{SIR2}
\pr\left(\zeta_{t}^{\{w,(\lambda,\mu)\}} \textup{ survives}\right)  \geq \pr\left(\xi_{t}^{\{w, \lambda\}} \textup{ survives}\right).
\end{equation}
\end{proposition}

The proof of this last statement relies on extending in our context an observation of Mollison \cite{M77} regarding a stochastic comparisson (coupling) between the standard spatial epidemic and a certain dependent directed percolation model, in which bonds are mutually independent if and only if they start out of different sites; this observation was subsequently developed in Kuulasmaa \cite{K82}. Our proof follows and extends the version of these arguments given in Durrett \cite{D88}, Chpt.\ 9. 

We remark next on an extension of monotonicity in the initial configuration property (\ref{rem11}) in juxtaposition with an important property of the contact process, known as additivity. This property provides that, if $\eta_{t}^{\eta}$ is the contact process, then 
\begin{equation}\label{addt}
\eta_{t}^{\eta_{1} \vee \eta_{2}} = \eta_{t}^{\eta_{1}} \vee \eta_{t}^{\eta_{2}},
\end{equation}
almost surely, where  $\eta_{1} \vee \eta_{2}$ denotes the configuration such that $\eta_{1} \vee \eta_{2}(x) = \max\{\eta_{1}(x), \eta_{2}(x)\}$, $x \in V$.  It may be easily shown that (\ref{addt}) is a stronger property than monotonicity in the initial configuration: $\eta_{1} \geq \eta_{2}  \Longrightarrow \mbox{ }\eta_{t}^{\eta_{1}} \geq_{st.} \eta_{t}^{\eta_{2}}$; see, for example, Corollary 1.3 in Chpt.\ II of Griffeath \cite{G79}. We point out that, although the analog of this property for the three state contact process is not known to hold, the following weaker conclusion can be deduced by the definition of stochastic domination and applying  (\ref{rem11}) twice.
\begin{remark}\textup{
Let $\mu \geq \lambda$. We have that: $\zeta_{t}^{\eta_{1} \vee \eta_{2}} \geq_{st.} \zeta_{t}^{\eta_{1}} \vee \zeta_{t}^{\eta_{2}}$, $\forall \eta_{1}, \eta_{2} \in \X$.}
\end{remark}

The next statement regards lack of monotonicity for the standard spatial epidemic, also known as the forest fire model.

\begin{proposition}\label{remmon}
Let $\xi_{t}^{\{A, \lambda\}} \equiv \zeta_{t}^{\{A, (\lambda,0)\}}$ and $\Xi_{t}^{A} = \{x: \xi_{t}^{\{A, \lambda\}}(x) = 1\}$. 
There is $G$ such that: 
\begin{equation}\label{rem1}
A \subset A' \hspace{5mm}  \centernot\Longrightarrow \mbox{ } \Xi_{t}^{\{A, \lambda\}} \subseteq_{st.} \Xi_{t}^{\{A', \lambda\}}, \mbox{for every } \lambda>0.  
\end{equation}
\begin{equation}\label{rem2}
\lambda < \lambda' \hspace{4mm} \centernot\Longrightarrow \mbox{ }\Xi_{t}^{\{A, \lambda\}} \subseteq_{st.} \Xi_{t}^{\{A, \lambda'\}}, \mbox{ for every } A \subset V. 
\end{equation}
\end{proposition} 

As far as we know, the statement in Proposition \ref{remmon} is not given explicitly  elsewhere in the literature. Intuitively, the result can be expected based on the following remarks, which we quote from $\mathsection$ 13.5 in Grimmett \cite{G99}. By adding an extra infective, one may subsequently infect a point which, during its removal period, prevents the infection from spreading 
further. A forest fire may be impeded by burning a pre-emptive firebreak. The proof of Proposition \ref{remmon} relies on constructing counterexamples on the connected graph with two vertices. We also point out here that a different approach for showing results alike Proposition \ref{remmon} would be through appropriate extensions of Theorem 3.2, Chpt.\ III in Liggett \cite{L85}, where the condition $\mu' \geq \mu$ is shown to be necessary and sufficient for attractiveness in the case of the contact process to hold.  The necessity part of this result relies on a simple argument relating the distribution of the process with its transition rates via its pregenarator, which however seems to not extend for the three state contact process. Finally, we note that the argument given in the proof of Proposition \ref{remmon} applies for $G_{1, +}$, giving in particular the following consequence mentioned above. 
\begin{remark}\label{conint}\textup{For $G \equiv G_{1, +} $ both conclusions of Proposition \ref{remmon} are valid. }
\end{remark}

\section{Proofs}

\begin{proof}[Proof of Theorem \ref{PIPmon}]   
We construct a coupled process $(\zeta_{t}'(x), \zeta_{t}(x))$ on $\X \times \X$ with the property that if $\eta'  \geq \eta$, then $\zeta_{t}' \geq \zeta_{t}$, for all $t \geq 0$, a.s., via prescribing all joint transitions-rates, i.e. flips. The existence of such a coupling implies (in fact it is equivalent to) the desired conclusion by a general result, see Theorem 2.4, Chpt.\ II, \cite{L85}. 
\begin{equation*}
(0, -1) \rightarrow \hspace{2mm}
\begin{cases}
(1,1) & \text{at rate } \lambda|\overrightarrow{yx} \in E : \zeta_{t}(y) = 1| \\
(1,-1) & \text{at rate } \mu' |\overrightarrow{yx} \in E : \zeta'_{t}(y) = 1| - \lambda |\overrightarrow{yx} \in E : \zeta_{t}(y) = 1|
\end{cases}
\end{equation*}
\begin{equation*}
(-1, -1) \rightarrow 
\begin{cases}
(1,1) & \text{at rate } \lambda|\overrightarrow{yx} \in E : \zeta_{t}(y) = 1| \\
(1,-1) & \text{at rate } \lambda' |\overrightarrow{yx} \in E : \zeta'_{t}(y) = 1| - \lambda |\overrightarrow{yx} \in E: \zeta_{t}(y) = 1|
\end{cases}
\end{equation*}

\begin{equation*}
(0, 0) \rightarrow \hspace{2mm}
\begin{cases}
(1,1) & \text{at rate } \mu|\overrightarrow{yx} \in E : \zeta_{t}(y) = 1| \\
(1,0) & \text{at rate } \mu' |\overrightarrow{yx} \in E : \zeta'_{t}(y) = 1| - \mu |\overrightarrow{yx} \in E : \zeta_{t}(y) = 1|
\end{cases}
\end{equation*}
Further, $(1, -1) \rightarrow (1,1)$ at rate $\lambda |\overrightarrow{yx} \in E : \zeta_{t}(y) = 1 |$, while 
$(1, -1) \rightarrow (0,-1)$ at rate 1. Also, $ (1, 0) \rightarrow (1,1)$ at rate $\mu |\overrightarrow{yx} \in E : \zeta_{t}(y) = 1 |$, while $ (1, 0) \rightarrow (0,0)$ at rate 1. Finally,  $(1, 1) \rightarrow (0,0)$ at rate 1. 

Note first that the inequality $\zeta_{t}'(x) \geq \zeta_{t}(x)$ is preserved by all flips prescribed above. Note also that the assumptions on the infection rates in (\ref{monimp}) are required for all of the flip rates of $(\zeta_{t}'(x), \zeta_{t}(x))$ prescribed above to be non-negative.  In addition, we observe that the flip rates define a valid coupling, since the marginals yield the correct transition rates for both $(\zeta_{t})$ and $(\zeta_{t}')$. To see this, one adds up the transition rates in order to check that the coordinates  $\zeta_{t}'(x)$ and  $\zeta_{t}(x)$  individually flip at the correct rates. For instance, from the first and third displays above, we have that $\zeta_{t}'(x): 0 \rightarrow 1$ at rate $\mu' |\overrightarrow{yx} \in E : \zeta'_{t}(y) = 1|$, whenever $\zeta_{t}(x)$ equals $-1$ or $0$ respectively. We have achieved the coupling and thus the proof is complete.
\end{proof}

\begin{proof}[Proof of Theorem \ref{onesided}] We first prove (\ref{eq:5}).  Let $I_n$ and $I'_n$ denote the sets of infected times for site $\{n\}$.
For each $n\geq 0$ these sets are denumerable unions of disjoint intervals.  
We are going to prove by induction on $n$ that there exists a representation
\begin{equation}\label{eq:3}
  I_n = \bigcup_{i\ge0} [a^{(n)}_i, b^{(n)}_i)
\end{equation}
where the index $i$ is used for the successive
infection intervals, and an associated (monotone increasing) function $\phi_n\colon{}I_n\to{}I'_n$ from $I_n$
into $I'_n$ such that, for all $i\ge0$,
\begin{align}
  \phi_n(a^{(n)}_i) & \le a^{(n)}_i \label{eq:1}\\
  \phi_n(t) & = \phi_n(a^{(n)}_i) + t - a^{(n)}_i,
  \qquad t \in [a^{(n)}_i,b^{(n)}_i)  \label{eq:2}\\
  \phi_n(a^{(n)}_{i+1})- \phi_n(a^{(n)}_i)
  & \ge  a^{(n)}_{i+1} -  a^{(n)}_{i}. \label{eq:4}
\end{align}

\begin{tikzpicture}[scale=7]
\draw[->, thick] (-0.1,0) -- (1.45,0);
\foreach \x/\xtext in {0.3/$a_{1}^{(n)}$,0.6/$b_{1}^{(n)}$,0.8/$a_{2}^{(n)}$,1.3/$b_{2}^{(n)}$}
    \draw[thick] (\x,0.5pt) -- (\x,-0.5pt) node[below] {\xtext};
\draw[[-), ultra thick, blue] (0.3,0) -- (0.6,0);
\draw[[-), ultra thick, green] (0.8,0) -- (1.3,0);
\draw (-0.25,0) node {$I_{n}$};
\end{tikzpicture}

\begin{tikzpicture}[scale=7]
\draw[->, thick] (-0.1,0) -- (1.45,0);
\foreach \x/\xtext in {0/$\phi_{n}(a_{1}^{(n)})$,0.3/$\phi_{n}(b_{1}^{(n)})$,0.6/$\phi_{n}(a_{2}^{(n)})$,1.1/$\phi_{n}(b_{2}^{(n)})$}
    \draw[thick] (\x,0.5pt) -- (\x,-0.5pt) node[below] {\xtext};
\draw[[-), ultra thick, blue] (0,0) -- (0.3,0);
\draw[[-), ultra thick, green] (0.6,0) -- (1.1,0);
\draw (-0.25,0) node {$I'_{n}$};
\end{tikzpicture}

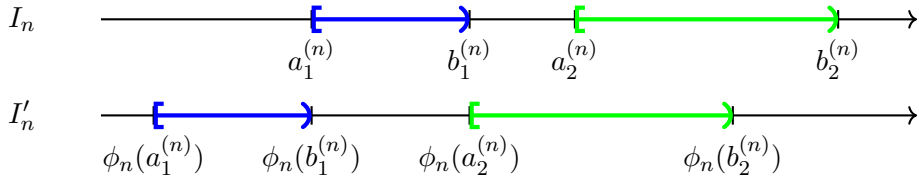
{\it\small\noindent Figure: Note that the image under $\phi_n$ of each of the intervals in the representation of $I_n$ in \eqref{eq:3} is an earlier interval in $I'_n$
of the same length, whereas the image intervals are disjoint and more widely spaced than the originals.}


We construct the coupled processes by induction on $n$.  Clearly 
the above representation in the case $n=0$ holds.
Therefore, we assume that we have such a representation for some $n$ 
and establish it for $n+1$. In particular this implies (\ref{eq:5}).  
Consider the $\zeta_{t}$ process first.  Given the set of
infected times $I_n$ at site $n$ we construct $I_{n+1}$ as follows.
On $I_n$ generate independent Poisson processes with rates $\lambda$,
corresponding to (potential) initial infections of site $n+1$, and $\mu$,
corresponding to secondary infections of site $n+1$.  Denote by
$\nu_0$ the time of the first point (in $I_n$) of the process at
rate $\lambda$, and by $\nu_1,\nu_2,\dots,\nu_k$ the times of
subsequent points (again in $I_n$) of the process at rate
$\mu$.  Thus $\nu_0\le\nu_1\le\dots\le\nu_k$ are the times at which
site $n+1$ is infected (if not already currently infected) from site
$n$.  Construct also an independent Poisson process with rate $1$ on
$\Rp$ defining the recovery events as the times at which site $n+1$ recovers.  Then the set of infected times $I_{n+1}$ at site
$n+1$ has the representation~\eqref{eq:3} with $a^{(n+1)}_i=\nu_i$,
$i=0,1,\dots,k$, and, for each such $i$, $b^{(n+1)}_i=\nu_i+d_i$,
where $\nu_i+d_i$ is the minimum of $\nu_{i+1}$ (where we define
$\nu_{k+1}=\infty$) and the time of the first recovery event after
$\nu_i$.

Now consider the $\zeta_{t}'$ process.  Given the set $I'_n$
of infected times at site $n$ we similarly construct $I'_{n+1}$ as
follows.  The independent Poisson processes of rates $\lambda'$ and
$\mu'$ on $I'_n$ (defining respectively the times of initial and secondary infections 
of site $n+1$) are given as
follows.  On the image $\phi_n(I_n)$ of $I_n$ in $I'_n$, these two
processes are given by using $\phi_n$ to map the points of the
corresponding Poisson processes (with rates $\lambda$ and $\mu$) on
$I_n$ which were used in the construction of the $\zeta_{t}$
process; in order to obtain the correct rates, these two processes are
then supplemented by the points of additional independent Poisson
processes of rates $\lambda'-\lambda$ and $\mu'-\mu$.  On
$I'_n\setminus{}\phi_n(I_n)$ we simply run additional independent Poisson
processes with rates $\lambda'$ and $\mu'$.  Denote by $\nu'_0$ the
time of the first point (in $I'_n$) of the process with rate
$\lambda'$, and, for $i=1,\dots,k$, define
$\nu'_i=\phi_n(\nu_i)\in{}I'_n$.  Thus $\nu'_0$ is the time of the
first infection (in the $\zeta_{t}'$ process) of site $n+1$,
while $\nu'_1,\dots,\nu'_k$ are a subset of the further times
at which site $n+1$ is infected (if not already currently infected)
from site $n$.  Since also, by construction, $\nu'_0\le\phi_n(\nu_0)$,
it follows from the properties of the function $\phi_n$ that
\begin{align}
  \nu'_i & \le \nu_i,
  \qquad i = 0,\dots,k,  \label{eq:6}\\
  \nu'_{i+1} - \nu'_i & \ge \nu_{i+1} - \nu_i,
  \qquad i = 0,\dots,k-1. \label{eq:7} 
\end{align}
The independent Poisson process with rate $1$ on $\Rp$ which defines
the times of the recovery events for infections of site $n+1$ in the
$\zeta_{t}'$ process is given as follows: consider the
corresponding Poisson process with rate $1$ used in the construction
of $I_{n+1}$.  The restriction of
this process to each of the intervals $[\nu_i,\nu_{i+1})$,
$i=0,1,\dots,k$, (and again with $\nu_{k+1}=\infty$) is mapped to the
interval $[\nu'_i,\nu'_i+\nu_{i+1}-\nu_i)$ in the obvious manner,
i.e.\ by adding $\nu'_i-\nu_i$ to each point (recall that, from
\eqref{eq:7}, these latter intervals are disjoint); outside the
intervals $[\nu'_i,\nu'_i+\nu_{i+1}-\nu_i)$ we place an independent
Poisson process with rate $1$.  The set $I'_{n+1}$ is now constructed
in the usual manner.  Note that, from the above construction, it
contains each of the intervals
$[\nu'_i,\nu'_i+d_i)\subseteq[\nu'_i,\nu'_i+\nu_{i+1}-\nu_i)$ for
$i=0,\dots,k$.  We can thus take the mapping $\phi_{n+1}$ to be given
by $\phi_{n+1}(a^{(n+1)}_i)=\phi_{n+1}(\nu_i)=\nu'_i$ for
$i=0,\dots,k$ and to be such that \eqref{eq:2} is satisfied with $n+1$
replacing $n$.  It follows from the above construction of $I_{n+1}$
and $I'_{n+1}$ that this indeed maps the former set into the latter.
Further it follows from \eqref{eq:6} and \eqref{eq:7} that
\eqref{eq:1} and \eqref{eq:4} are similarly satisfied with $n+1$
replacing $n$. This proves (\ref{eq:5}). The proof of (\ref{extra}) is omitted as it follows from (\ref{eq:5}) by an instance of the argument given in the the last paragraph of the proof of Proposition \ref{propff} for deriving (\ref{SIR2}) from (\ref{SIR1}) there.  This completes the proof.  

\end{proof}

\begin{proof}[Proof of Proposition \ref{propff}] 
We first show that we may define the two processes on a common probability space, such that the following holds: 
\begin{equation}\label{SIR1}
\left\{\xi_{t}^{\{w, \lambda\}}(v) = 1, \mbox{ for some } t\geq0  \right\} \subseteq \left\{\zeta_{t}^{\{w,(\lambda,\mu)\}}(v) = 1, \mbox{ for some } t\geq0  \right\}, 
\end{equation}
almost surely. To do so, for all $u\in V$ let $(T_{n}^{u})_{n\geq1}$ be exponential 1 r.v.'s; further for all $u, v$ such that $\overrightarrow{uv}$, let $(Y_{n}^{(u,v)})_{n\geq1}$ be exponential $\lambda$  r.v.'s and $(N_{n}^{(u,v)})_{n\geq1}$ be Poisson processes at rate $\mu$. All random elements introduced are independent and $\pr$ below denotes the corresponding probability measure. To describe the construction below, let $\tau_{k,n}^{(u,v)}, k\geq1$, be the times of events of $N_{n}^{(u,v)}$ within the time interval $[0,T_{n}^{u})$ and let also $X_{n}^{(u,v)}, n\geq1$, be such that $X_{n}^{(u,v)} = Y_{n}^{(u,v)} $ if $Y_{n}^{(u,v)} < T_{n}^{u}$ and $X_{n}^{(u,v)} := \infty$ otherwise.

We now construct $\zeta_{t}^{\{w,(\lambda,\mu)\}}$ on $G$ as follows. Suppose that site $u$ gets infected at time $t$ for the $n$-th time, $n\geq1$, then: (i) at time $ t + T_{n}^{u}$ a recovery occurs at site $u$, (ii) at time $t+X_{n}^{(u,v)}$ an initial infection of $v$ occurs if immediately prior to that time site $v$ is at a susceptible and never infected state, and, (iii) at each time $t+\tau_{k,n}^{(u,v)}, k\geq1$, a secondary infection occurs at site $v$ if immediately prior to that time site $v$ is at a susceptible and previously infected state.

\[\begin{tikzcd}[%
row sep = 5mm, column sep = 5mm, 
cells={nodes={circle, draw}}, 
every arrow/.append style={shorten <= 1mm, shorten >= 1mm}] 
1  \arrow[shift left]{d} \arrow[shift left]{r} & 2 \arrow[shift left]{d} \arrow[shift left]{l} \arrow[shift left]{r} & 3 \arrow[shift left]{d} \arrow[shift left]{l} & 1  \arrow[shift left]{d} \arrow[shift left]{r} & 2 \arrow[shift left]{d} \arrow[shift left]{r} & 3 \arrow[shift left]{l} \\
4 \arrow[shift left]{u}  \arrow[shift left]{d} \arrow[shift left]{r} & 5 \arrow[shift left]{l} \arrow[shift left]{u} \arrow[shift left]{d} \arrow[shift left]{r} &  6 \arrow[shift left]{l} \arrow[shift left]{u} \arrow[shift left]{d} &  4  \arrow[shift left]{r} & 5\arrow[shift left]{d} &  6 \arrow[shift left]{u}  \\
7 \arrow[shift left]{u} \arrow[shift left]{r} & 8 \arrow[shift left]{l} \arrow[shift left]{u} \arrow[shift left]{r} &  9 \arrow[shift left]{l}  \arrow[shift left]{u} &  7 \arrow[shift left]{u}  & 8 \arrow[shift left]{l} \arrow[shift left]{r}  &  9  
\end{tikzcd}\]
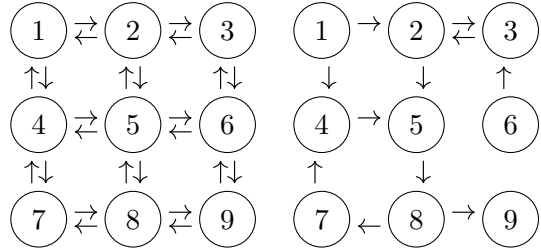
{\it\small Figure: An example of a digraph $G$ (left); a realization of the associated (dependent) directed percolation random graph $\Gamma$ (right), in which, for instance, $\mathcal{X}_{2} = \{3,5\}$ and $1 \mbox{ }_{\overrightarrow{(\mathcal{X}_{u}, u\in  V)}}  \mbox{ }9 $.}

Let also $\mathcal{X}_{u} = \{v: \overrightarrow{uv} \mbox{ and } X_{1}^{(u,v)}<\infty\}$, $u\in V$. Let  $\Gamma$ denote the subgraph of $G$ induced by retaining edges from $u$ to $v$ if and only if $v \in \mathcal{X}_{u}$, for all $u,v \in  V$. Let further $u \mbox{ }_{\overrightarrow{(\mathcal{X}_{u}, u\in  V)}}  \mbox{ }v $ denote the existence of a directed path from $u$ to $v$ in $\Gamma$. It follows by Lemma 1 in \cite{D88}, Chpt.\ 9, that by the construction of $\zeta_{t}^{\{w,(\lambda,0)\}}$, we have that 
\begin{equation*}\label{Czeta2}
\{w \mbox{ }_{\overrightarrow{(\mathcal{X}_{u}, u\in  V)}}  \mbox{ } v \}  = \left\{\zeta_{t}^{\{w,(\lambda,0)\}}(v) = 1  \mbox{ for some } t\geq0  \right\}. 
\end{equation*}
Analogously for $\zeta_{t}^{\{w,(\lambda,\mu)\}}$, we also have that 
\begin{equation*}\label{Czeta}
\{w \mbox{ }_{\overrightarrow{(\mathcal{X}_{u}, u\in  V)}}  \mbox{ } v \}  \subseteq \left\{\zeta_{t}^{\{w,(\lambda,\mu)\}}(v) = 1  \mbox{ for some } t\geq0  \right\}, 
\end{equation*}
for all $v\in  V$. The proof of (\ref{SIR1}) is then complete by combining the last two displays above. 

We finally derive (\ref{SIR2}) from (\ref{SIR1}). We consider $\zeta_{t}^{\{w,(\lambda,\mu)\}}$ and for $v \in  V$, we let 
\[
A_{v} = \{\zeta_{t}^{\{w,(\lambda,\mu)\}}(v) =1\mbox{ for some } t\geq0\}. 
\]
From the first part, it suffices to show that
\begin{equation*}\label{onsurvinfclust}
\pr\left(\sum_{v \in  V} 1(A_{v}) = \infty \right) = \pr\big(\zeta_{t}^{\{w,(\lambda,\mu)\}}\mbox{ survives}\big),
\end{equation*}
where $1(\cdot)$ denotes the indicator function. To prove the equality in the last display above, we let $B_{M}$ denote the event  $\big\{ \sum_{v \in  V} 1(A_{v}) \leq M \big\}$. Either because a finite state-space irreducible Markov chain with a single absorbing state is eventually absorbed, or by Lemma VII.4.1 in \cite{S99}, we have that, for all fixed integer $M\in [1, \infty)$, 
\[
\pr\big(B_{M}, \zeta_{t}^{\{w,(\lambda,\mu)\}}\mbox{ survives}\big) =0,
\] 
and thus, $\pr\Big(\bigcup_{M\geq1} B_{M},  \zeta_{t}^{\{w,(\lambda,\mu)\}}\mbox{ survives}\Big) = 0,$ which completes the proof. 
\end{proof}

\begin{proof}[Proof of Proposition \ref{remmon}] Let $G$ be the connected graph with $V = \{u,v\}$. We will show  that: $(i)$ for all $\lambda>1$, a coupling of $\zeta_{t}^{\{u, (\lambda,0)\}}$ and $\zeta_{t}^{\{V, (\lambda,0)\}}$ on $G$, 
such that $\zeta_{t}^{\{u, (\lambda,0)\}} \leq \zeta_{t}^{\{V, (\lambda,0)\}}, \forall \mbox{ } t \geq0$, cannot be constructed; and, further that:  $(ii)$ for all $\lambda, \lambda'$, if $\lambda < \lambda'<1$ then  a coupling of $\zeta_{t}^{\{u, (\lambda,0)\}}$ and $\zeta_{t}^{\{u, (\lambda',0)\}}$ on $G$ such that $\zeta_{t}^{\{u, (\lambda,0)\}} \leq \zeta_{t}^{\{u, (\lambda',0)\}}, \forall \mbox{ } t \geq0$, cannot be constructed. This suffices since $(i)$ and $(ii)$  imply $(\ref{rem1})$ and $(\ref{rem2})$ respectively.

Let $T_{u},T_{v}$ be exponential 1 r.v.'s; let also $X_{u,v}$ be an exponential $\lambda$ \mbox{r.v.}, and $f_{X_{u,v}}$ be its probability density function. All r.v.'s introduced are independent of each other and defined on some probability space with probability measure $\pr$. We have that, for any $t\geq0$,  
\begin{eqnarray}\label{calc}
\pr\left(\zeta_{t}^{\{u, (\lambda,0)\}} = (1,1) \right) & = &  \pr(T_{u}>t) \int_{0}^{t} f_{X_{u,v}}(s)\pr(T_{v}> t-s)\, ds \nonumber\\
& = & e^{-2t}  \int_{0}^{t} \lambda e^{s (1- \lambda)} \, ds \nonumber \\
& =& e^{-2t}  \frac{\lambda}{\lambda-1} (1-e^{-t(\lambda-1)}).
\end{eqnarray}

By $(\ref{calc})$ then we have: (a) for all $\lambda >1$ we can choose $t$ sufficiently large, i.e. $\displaystyle{ t > \frac{\log{\lambda}}{\lambda-1}}$, such that $\pr\left(\zeta_{t}^{\{u, (\lambda,0)\}}=(1,1)\right) > e^{-2t} =  \pr\left(\zeta_{t}^{\{V, (\lambda,0)\}} = (1,1) \right)$; and further that (b) for all $\lambda <1$, $\pr\left(\zeta_{t}^{\{u, (\lambda,0)\}} = (1,1)\right)$ is not an increasing function of $\lambda$. From Theorem B9 in \cite{L99}, $(i)$ and $(ii)$ follow by (a) and (b) respectively. This completes the proof.  

\end{proof} 

\section*{Acknowledgments}
Thanks to two anonymous referees for useful comments and for a careful reading of the manuscript. This work was partly supported during non-overlapping periods of time: by Heriot-Watt University, by CONICET, by part of FAPESP, project, titled Research, Innovation and Dissemination Center for Neuromath (grant 2013/ 07699-0), FAPESP grant 2016/03988-5, and by PNPD/CAPES grant.

\end{document}